\documentclass[11pt]{article}
\usepackage{amssymb,amsmath,amsthm}
\newcommand{\R}{\mathbb{R}}
\newcommand{\Z}{\mathbb{Z}}
\newcommand{\N}{\mathbb{N}}
\def\qed{\hfill $\Box$ \smallskip}
\def\x#1{{\rm (\ref{#1})}}

\begin{document}
\title{\LARGE\bf{ New super-quadratic conditions for asymptotically periodic Schr\"odinger equation}
 \footnote{This work is partially supported by the NNSF (No: 11171351) and the SRFDP(No: 20120162110021) of China}}
\date{}
 \author{ X. H. Tang\\
        {\small School of Mathematics and Statistics,}\\
        {\small Central South University,}\\
        {\small Changsha, Hunan 410083, P.R.China }\\
        {\small E-mail: tangxh@mail.csu.edu.cn}}
\maketitle
\begin{center}
\begin{minipage}{13cm}
\par
\small  {\bf Abstract:} This paper is dedicated to studying the semilinear Schr\"odinger equation
 $$
   \left\{
   \begin{array}{ll}
    -\triangle u+V(x)u=f(x, u), \ \ \ \  x\in {\R}^{N},\\
    u\in H^{1}({\R}^{N}),
   \end{array}
    \right.
 $$
 where $f$ is a superlinear, subcritical nonlinearity. It focuses on the case where
 $V(x)=V_0(x)+V_1(x)$, $V_0\in C(\R^N)$, $V_0(x)$ is 1-periodic in each of $x_1, x_2, \ldots, x_N$ and
 $\sup[\sigma(-\triangle +V_0)\cap (-\infty, 0)]<0<\inf[\sigma(-\triangle +V_0)\cap (0, \infty)]$,
 $V_1\in C(\R^N)$ and $\lim_{|x|\to\infty}V_1(x)=0$. A new super-quadratic condition is obtained,
 which is weaker than some well known results.

 \vskip2mm
 \par
 {\bf Keywords: }  Schr\"odinger equation; Superlinear; Asymptotically periodic;
 Ground state solutions of Nehari-Pankov type.

 \vskip2mm
 \par
 {\bf 2000 Mathematics Subject Classification.}  35J20; 35J60
\end{minipage}
\end{center}

 {\section{Introduction}}
 \setcounter{equation}{0}

 \par
   Consider the following semilinear Schr\"odinger equation
 \begin{equation} \label{SP}
   \left\{
   \begin{array}{ll}
    -\triangle u+V(x)u=f(x, u), \ \ \ \  x\in {\R}^{N},\\
    u\in H^{1}(\R^N),
   \end{array}
    \right.
 \end{equation}
 where $V : {\R}^{N} \rightarrow {\R}$ and $f: {\R}^N\times {\R} \rightarrow {\R}$ are asymptotically periodic
 in $x$, moreover $f$ is super linear as $|u| \rightarrow \infty$.

  \par
   When $V(x)$ and $f(x, u)$ are periodic in $x$, and satisfy the following basic assumptions, respectively:

 \vskip2mm
 \noindent
 \begin{itemize}
 \item[(V)] $V\in C(\R^N)$, $V(x)$ is 1-periodic in each of $x_1, x_2, \ldots, x_N$ and
 \begin{equation}\label{sp}
  \sup[\sigma(-\triangle +V)\cap (-\infty, 0)]<0<\bar{\Lambda}:=\inf[\sigma(-\triangle +V)\cap (0, \infty)];
 \end{equation}

 \item[(F1)] $f\in C(\R^N\times \R)$, and there exist constants $p\in (2, 2^*)$ and $C_0>0$ such that
 $$
   |f(x, t)|\le C_0\left(1+|t|^{p-1}\right), \ \ \ \ \forall \ (x, t)\in \R^N\times \R;
 $$

 \item[(F2)] $f(x, t)=o(|t|)$, as $|t|\to 0$, uniformly in $x\in \R^{N}$, and $F(x, t):=\int_{0}^{t}f(x, s)\mathrm{d}s\ge 0$;

 \item[(F3)] $f(x, t)$ is 1-periodic in each of $x_1, x_2, \ldots, x_N$;
 \end{itemize}

 \vskip2mm
 \noindent
 the existence of a nontrivial solution for \x{SP} has been widely investigated in literature, for example, see
 \cite{BW, CZR, KS, LS, Ra, TW, WZ} and references cited therein.
 In these papers, a classical existence condition is (AR) which is due to Ambrosetti and Rabinowitz \cite{AR}:

 \vskip2mm
 \noindent
 \begin{itemize}
 \item[(AR)] there exists a $\mu>2$ such that
 $$
   0<\mu F(x, t)\le tf(x, t), \ \ \ \ \forall \ (x, t)\in \R^N\times (\R\setminus \{0\}).
 $$
 \end{itemize}

 \vskip2mm
 \par
  (AR) is a very convenient hypothesis since it readily achieves mountain pass geometry as well
 as satisfaction of the Palais-Smale condition. However, it is a severe restriction, since it strictly controls the
 growth of $f(x, t)$ as $|t|\to \infty$. In recent years, there are some papers devoted to replace (AR) with
 weaker conditions. For example, Liu and Wang \cite{LW} first introduced a more natural super-quadratic
 condition:

 \vskip2mm
 \noindent
 \begin{itemize}
 \item[(SQ)] $\lim_{|t|\to \infty}\frac{|F(x, t)|}{|t|^2}=\infty$, uniformly in $x\in \R^{N}$.
 \end{itemize}

 \vskip2mm
 %\noindent
 Subsequently, it has been commonly used in many recent papers, see \cite{DL, DL1, DS, LS, LWZ, Pa, Ta, Ya}. However, to some extend,
 the condition(SQ) also has its own limitation, for it is not sufficient to guarantee \x{SP} has a nontrivial solution.

 \par
   Later, in 2006, Ding and Lee \cite{DL} gave a more mild existence condition:

 \vskip2mm
 \noindent
 \begin{itemize}
 \item[(DL)]  $\mathcal{F}(x, t):=\frac{1}{2}tf(x, t)-F(x, t) > 0$ if $t\ne 0$, and there exist $c_0>0$, $r_0>0$ and $\kappa>$
 $\max\{1, N/2\}$ such that
 $$
   |f(x, t)|^{\kappa}\le c_0\mathcal{F}(x, t) |t|^{\kappa},  \ \ \ \ \forall \ (x, t)\in \R^N\times \R, \ \ |t|\ge r_0.
 $$
 \end{itemize}

 \vskip2mm
 \noindent
 Under the assumption (F1), the condition (DL) greatly weaken (AR). Soon after, it was generalized in various directions and applied to more general
 equations or systems by numerous of authors, see e.g. \cite{BD, Ta, Ta2, YCD, ZCZ}.

 \par
   In paper \cite{Sz}, Szulkin and Weth developed an ingenious approach to find the ground state solutions
 for problem \x{SP}. They demonstrated that (SQ) together with the following Nehari type assumption (Ne)
 implies \x{SP} has a ground state solution.

 \vskip2mm
 \noindent
 \begin{itemize}
 \item[(Ne)] $t\mapsto f(x, t)/|t|$ is strictly increasing on $(-\infty, 0)\cup (0, \infty)$.
 \end{itemize}

 \vskip2mm
 \noindent
 Based on Szulkin and Weth \cite{Sz}, Liu \cite {Liu} showed that \x{SP} has a nontrivial solution by using the
 following weak version (WN) instead of (Ne):

 \vskip2mm
 \noindent
 \begin{itemize}
 \item[(WN)] $t\mapsto f(x, t)/|t|$ is non-decreasing on $(-\infty, 0)\cup (0, \infty)$.
 \end{itemize}

 \vskip2mm
 %\noindent
 In a very recent paper \cite{Ta1}, Tang introduced new super--quadratic conditions as follows:

 \vskip2mm
 \noindent
 \begin{itemize}
 \item[(WS)] $\lim_{|t|\to \infty}\frac{|F(x, t)|}{|t|^2}=\infty, \ \ a.e. \ x\in \R^N$;
 \end{itemize}

 \begin{itemize}
 \item[(Ta)] there exists a $\theta_0\in (0, 1)$ such that
 $$
   \frac{1-\theta^2}{2}tf(x, t) \ge \int_{\theta t}^{t}f(x, s)\mathrm{d}s=F(x, t)-F(x, \theta t),
       \ \ \ \ \forall \ \theta\in [0, \theta_0], \ \ (x, t)\in \R^N\times \R.
 $$
 \end{itemize}

 \vskip2mm
 \noindent
 Clearly, (WS) is slightly weaker than (SQ). Besides, (Ta) improves (AR), (WN) and a weak version of (AR) (see \cite{Ta1}):

 \vskip2mm
 \noindent
 \begin{itemize}
 \item[(WAR)] there exists a $\mu>2$ such that
 $$
   0\le \mu F(x, t)\le tf(x, t), \ \ \ \ \forall \ (x, t)\in \R^N\times \R.
 $$
 \end{itemize}

 \vskip2mm
 \par
   Motivated by the aforementioned works, in the periodic case, we first weaken (DL) to the following condition
 (i.e. $\mathcal{F}(x, t) > 0, t\ne 0$ to $\mathcal{F}(x, t) \ge 0$):

 \vskip4mm
 \noindent
 \begin{itemize}
 \item[(F4)] $\mathcal{F}(x, t) \ge 0$, and there exist $c_0>0$, $\delta_0\in (0, \bar{\Lambda})$ and $\kappa>\max\{1, N/2\}$
 such that
 $$
   \frac{f(x, t)}{t}\ge \bar{\Lambda}-\delta_0 \ \Rightarrow \ \left[\frac{f(x, t)}{t}\right]^{\kappa}\le c_0\mathcal{F}(x, t).
 $$
 \end{itemize}

 \vskip2mm
 \par
   Clearly, (WAR) and (DL) yield (F4). What we do notice, though, is that we cann't verify that (WN) implies (F4),
 it is very difficult to find a function $f$ which satisfies (F2) and (WN) but not (F4). Before presenting our first
 result, we give two nonlinear examples to illustrate the assumption (F4).

 \vskip4mm
 \par\noindent
 {\bf Example 1.1.}\ \ Let $F(x, t)=t^2\ln [1+t^2\sin^2(2\pi x_1)]$. Then
 $$
   f(x, t)=2t\ln \left[1+t^2\sin^2(2\pi x_1)\right]+\frac{2t^3\sin^2(2\pi x_1)}{1+t^2\sin^2(2\pi x_1)}, \ \ \ \
   \mathcal{F}(x, t)=\frac{t^4\sin^2(2\pi x_1)}{1+t^2\sin^2(2\pi x_1)}\ge 0.
 $$
 It is easy to see that $f$ does not satisfy (AR), (SQ), (WAR) and (DL), but it satisfies (WS) and (F4)
 with $\kappa>\max\{1, N/2\}$.

 \vskip4mm
 \par\noindent
 {\bf Example 1.2.}\ \ Let $N\le 4$ and
 $$
   F(x, t)=a\left(|t|^{13/4}-\frac{5}{2}|t|^{11/4}+\frac{45}{16}|t|^{9/4}\right), \ \ \ \ a>0.
 $$
 Then
 $$
   f(x, t)=a\left(\frac{13}{4}|t|^{5/4}-\frac{55}{8}|t|^{3/4}+\frac{405}{64}|t|^{1/4}\right)t,
 $$
 $$
   \mathcal{F}(x, t)=\frac{5}{8}a|t|^{9/4}\left(\sqrt{|t|}-\frac{3}{4}\right)^2\ge 0.
 $$
 Similarly, $f$ does not satisfy (AR), (WN), (WAR), (DL) and (Ta), but it satisfies (SQ) and (F4)
 with $\kappa=12/5$ if $a\in (0, 64\bar{\Lambda}/405)$.

 \vskip4mm
 \par
   We are now in a position  to state the first result of this paper.

 \vskip4mm
 \par\noindent
 {\bf Theorem 1.3.}\ \ {\it Assume that $V$ and $f$ satisfy} (V), (F1), (F2), (F3), (F4) {\it and}
 (WS). {\it Then problem \x{SP} has a nontrivial solution.}

 \vskip4mm
 \par
   When $V(x)$ is positive and asymptotically periodic, there are considerably fewer results;
 here we mention \cite{AD, LT, ZXZ}. In this case, the spectrum $\sigma (-\triangle +V)\subset (0, \infty)$.
 Comparing with appropriate solutions of a periodic problem
associated with  \x{SP},  a nontrivial solution was found by using a version of
the mountain pass theorem.
 \par
   When $V(x)$ is periodic and sign-changing, while $f(x, u)$ asymptotically periodic in $x$, there seems to be only
 one result \cite{LS}. Let $\Phi_0$ and $\Phi$ denote the energy functionals associated with problem \x{SP} with
 periodic and asymptotically periodic nonlinearity $f$, respectively. According to a generalized linking theorem for the strongly indefinite functionals, by comparing with (C)c-sequences of $\Phi_0$
 and $\Phi$, Li and Szuklin \cite{LS}
 proved that problem \x{SP} has a nontrivial solution if $V$ and $f$ satisfy assumptions (V) and (F1), (F2),
 (AR) and the following asymptotically periodic condition:

 \vskip2mm
 \noindent
 \begin{itemize}
 \item[(F5)] $f(x, t)=f_0(x, t)+f_1(x, t)$, $\partial_{t}f_0, f_1\in C(\R^N\times \R)$, $f_0(x, t)$ is 1-periodic in each
 of $x_1, x_2,$ $ \ldots, x_N$; and $f_0$ and $f_1$ satisfy that
 $$
   0 < tf_0(x, t)<t^2\partial_{t}f_0(x, t),  \ \ \ \ \forall \ (x, t)\in \R^N\times (\R\setminus \{0\}),
 $$
 $$
   0<F_0(x, t):=\int_{0}^{t}f_0(x, s)\mathrm{d}s\le \frac{1}{\mu}tf_0(x, t), \ \ \ \ \forall \ (x, t)\in \R^N\times (\R\setminus \{0\}),
 $$
 $$
   F_1(x, t):=\int_{0}^{t}f_1(x, s)\mathrm{d}s > 0, \ \ \ \ \forall \ (x, t)\in \R^N\times (\R\setminus \{0\})
 $$
 and
 $$
   |f_1(x, t)|\le a(x)\left(|t|+|t|^{p-1}\right), \ \ \ \ \forall \ (x, t)\in \R^N\times \R,
 $$
 where $\mu>2$ is the same as in (AR), $a\in C({\R}^{N})$ with $\lim_{|x|\to\infty}a(x)=0$.
 \end{itemize}

 \vskip2mm
 \par
    In the following, we will point out that the assumption that $f_0(x,t)$ is differentiable in $t$ and
 $0<tf_0(x, t)<t^2\partial_{t}f_0(x, t)$ in (F5) (which implies that $t\mapsto
 f_0(x, t)/|t|$ is strictly increasing on $(-\infty, 0)\cup(0, \infty)$) is very crucial in
 Li and Szulkin \cite{LS}.

 \par
   If $V(x)$ is both asymptotically periodic and sign-changing, the operator $-\triangle +V$ loses the $\Z^N$-translation
 invariance. For this reason, many effective methods for periodic problems cannot be applied to asymptotically periodic
 ones,  and they all failed. To the best of our knowledge, there are no existence results for \x{SP} when $V(x)$ is asymptotically periodic
 and sign-changing. Motivated by the works \cite{DL, LS, LT, Ta3, Ta4, Ta5, ZXZ}, we shall find new tricks
 to overcome the difficulties caused by the dropping of periodicity of $V(x)$.

 \vskip2mm
 \par
   Before presenting our second theorem, we make the following assumptions instead of (V) and (F5), respectively.

 \vskip2mm
 \noindent
 \begin{itemize}
 \item[(V1)] $V(x)=V_0(x)+V_1(x)$, $V_0\in C(\R^N)\cap L^{\infty}(\R^N)$ and
 \begin{equation}\label{sp0}
  \sup[\sigma(-\triangle +V_0)\cap (-\infty, 0)]<0< \bar{\Lambda}:=\inf[\sigma(-\triangle +V_0)\cap (0, \infty)],
 \end{equation}
 $V_1\in C(\R^N)$ and $\lim_{|x|\to\infty}V_1(x)=0$;

 \item[(V2)] $V_0(x)$ is 1-periodic in each of $x_1, x_2, \ldots, x_N$, and
 $$
   0\le -V_1(x)\le \sup_{\R^N}[-V_1(x)]<\bar{\Lambda}, \ \ \ \ \forall \ x\in \R^N;
 $$

 \item[(F5$'$)] $f(x, t)=f_0(x, t)+f_1(x, t)$, $f_0\in C(\R^N\times \R)$, $f(x, t)$ is 1-periodic in each of
 $x_1, x_2, \ldots, x_N$, $f_0(x, t)=o(|t|)$, as $t\to 0$, uniformly in $x\in \R^N$, $t\mapsto f_0(x, t)/|t|$
 is non-decreasing on $(-\infty, 0)\cup (0, \infty)$,
 $\lim_{|t|\to \infty}\frac{|F_0(x, t)|}{|t|^2}=\infty, \ \ a.e. \ x\in \R^N$; $f_1\in C(\R^N\times \R)$ satisfies that
 $$
   -V_1(x)t^2+F_1(x, t) > 0, \ \ \ \ \forall \ (x, t)\in B_{1+\sqrt{N}}(0)\times (\R\setminus \{0\})
 $$
 and
 $$
   F_1(x, t)\ge 0, \ \ \ \ |f_1(x, t)|\le a(x)\left(|t|+|t|^{p-1}\right), \ \ \ \ \forall \ (x, t)\in \R^N\times \R,
 $$
 where $a\in C({\R}^{N})$ with $\lim_{|x|\to\infty}a(x)=0$.
 \end{itemize}

 \vskip4mm
 \par\noindent
 {\bf Remark 1.4.}\ \ {\it Compare } (F5) {\it with} (F5$'$), {\it it is dropped that $f_0(x,t)$ is
 differentiable in $t$, moreover, the condition that $t\mapsto f_0(x, t)/|t|$ is strictly increasing on
 $(-\infty, 0)\cup(0, \infty)$ is weaken to non-decreasing on $(-\infty, 0)\cup(0, \infty)$; and that
 $f_0$ satisfies} (AR) {\it is also weaken to} (WS).

 \vskip4mm
 \par
   We are now in a position to state the second result of this paper.

 \vskip4mm
 \par\noindent
 {\bf Theorem 1.5.}\ \ {\it Assume that} (V1), (V2), (F1), (F2), (F4) {\it and} (F5$'$) {\it are satisfied. Then problem
 \x{SP} has a nontrivial solution.}

 \vskip4mm
 \par
   The remainder of this paper is organized as follows. In Section 2, some preliminary results are presented.
 The proofs of Theorems 1.3 and 1.5 are given in Section 3 and Section 4, respectively.

 \vskip10mm

 {\section{Preliminaries}}
 \setcounter{equation}{0}

 \vskip4mm
 \par
   Let $X$ be a real Hilbert space with $X=X^{-}\oplus X^{+}$ and $X^{-}\bot\  X^{+}$.
 For a functional $\varphi\in C^{1}(X, \R)$, $\varphi$ is said to be weakly sequentially
 lower semi-continuous if for any $u_n\rightharpoonup u$ in $X$ one has $\varphi(u)\le \liminf_{n\to\infty}\varphi(u_n)$,
 and $\varphi'$ is said to be weakly sequentially continuous if $\lim_{n\to\infty}\langle\varphi'(u_n), v\rangle=
 \langle\varphi'(u), v\rangle$ for each $v\in X$.

 \vskip4mm
 \par\noindent
 {\bf Lemma 2.1.} (\cite {KS, LS})\ \ {\it Let $(X, \|\cdot\|)$ be a real Hilbert space with $X=X^{-}\oplus X^{+}$ and $X^{-}\bot\  X^{+}$,
 and let $\varphi\in C^{1}(X, \R)$ of the form
 $$
   \varphi(u)=\frac{1}{2}\left(\|u^{+}\|-\|u^{-}\|\right)-\psi(u), \ \ \ \ u=u^{-}+u^{+}\in X^{-}\oplus X^{+}.
 $$
 Suppose that the following assumptions are satisfied:}

 \vskip2mm
 \noindent
 \begin{itemize}
 \item[(KS1)] {\it $\psi\in C^{1}(X, \R)$ is bounded from below and weakly sequentially lower semi-continuous;}

 \item[(KS2)] {\it $\psi'$ is weakly sequentially continuous;}

 \item[(KS3)] {\it there exist $r>\rho>0$ and $e\in X^{+}$ with $\|e\|=1$ such that
 $$
   \kappa:=\inf\varphi(S^{+}_{\rho}) > \sup \varphi(\partial Q),
 $$
 where }
 $$
     S^{+}_{\rho}=\left\{u\in X^{+} : \|u\|=\rho\right\},  \ \ \ \ Q=\left\{v+se : v\in X^{-},\ s\ge 0,\ \|v+se\|\le r\right\}.
 $$
 \end{itemize}

 \vskip2mm
 \noindent
 {\it Then there exist a constant $c\in [\kappa, \sup \varphi(Q)]$ and a sequence $\{u_n\}\subset X$ satisfying}
 \begin{eqnarray*}
   \varphi(u_n)\rightarrow c, \ \ \ \ \|\varphi'(u_n)\|(1+\|u_n\|)\rightarrow 0.
 \end{eqnarray*}

 \vskip4mm
 \par
   Let $\mathcal{A}_0=-\triangle +V_0$. Then $\mathcal{A}_0$ is self-adjoint in $L^2(\R^N)$ with domain
 $\mathfrak{D}(\mathcal{A}_0)=H^2(\R^N)$ (see \cite[Theorem 4.26]{EK}). Let $\{\mathcal{E}(\lambda):
 -\infty \le \lambda \le +\infty\}$ and $|\mathcal{A}_0|$ be the spectral family and the absolute value of $\mathcal{A}_0$,
 respectively, and $|\mathcal{A}_0|^{1/2}$ be the square root of $|\mathcal{A}_0|$. Set $\mathcal{U}=id-\mathcal{E}(0)
 -\mathcal{E}(0-)$. Then $\mathcal{U}$ commutes with $\mathcal{A}_0$, $|\mathcal{A}_0|$ and $|\mathcal{A}_0|^{1/2}$, and
 $\mathcal{A}_0 = \mathcal{U}|\mathcal{A}_0|$ is the polar decomposition of $\mathcal{A}_0$ (see \cite[Theorem IV 3.3]{EE}).
 Let
 \begin{equation}\label{E}
   E=\mathfrak{D}(|\mathcal{A}_0|^{1/2}), \ \ \ \ E^{-}=\mathcal{E}(0)E, \ \ \ \ E^{+}=[id-\mathcal{E}(0)]E.
 \end{equation}
 For any $u\in E$, it is easy to see that $u=u^{-}+u^{+}$, where
 \begin{equation}\label{u-}
   u^{-}:=\mathcal{E}(0)u\in E^{-},  \ \ \ \ u^{+}:=[id-\mathcal{E}(0)]u\in E^{+}
 \end{equation}
 and
 \begin{equation}\label{Au}
   \mathcal{A}_0u^{-}=-|\mathcal{A}_0|u^{-}, \ \ \ \ \mathcal{A}_0u^{+}=|\mathcal{A}_0|u^{+},
       \ \ \ \ \forall \ u\in E\cap \mathfrak{D}(\mathcal{A}_0).
 \end{equation}

 \par
    Define an inner product
 \begin{equation} \label{IP}
    (u, v)=\left(|\mathcal{A}_0|^{1/2}u, |\mathcal{A}_0|^{1/2}v\right)_{L^2}, \ \ \ \ u, v\in E
 \end{equation}
 and the corresponding norm
 \begin{equation} \label{No0}
    \|u\|=\left\||\mathcal{A}_0|^{1/2}u\right\|_{2}, \ \ \ \ u\in E,
 \end{equation}
 where $(\cdot, \cdot)_{L^2}$ denotes the inner product of $L^2(\R^N)$, $\|\cdot\|_s$ denotes the norm of $L^s(\R^N)$.
 By (V1), $E=H^1(\R^N)$ with equivalent norms. Therefore, $E$ embeds continuously in $L^s(\R^N)$ for all
 $2\le s\le 2^*$. In addition, one has the decomposition $ E=E^{-}\oplus E^{+}$ orthogonal with respect to
 both $(\cdot, \cdot)_{L^2}$ and $(\cdot, \cdot)$.

 \par
   Under assumptions (V1), (F1) and (F2), the solutions of problem \x{SP} are critical points of the functional
 \begin{equation}\label{Ph}
   \Phi(u)=\frac{1}{2}\int_{{\R}^N}\left(|\nabla u|^2+V(x)u^2\right)\mathrm{d}x-\int_{{\R}^N}F(x, u)\mathrm{d}x,
      \ \ \ \ \forall \ u\in E,
 \end{equation}
 $\Phi$ is of class $C^{1}(E, \R)$, and
 \begin{equation}\label{Phd}
   \langle \Phi'(u), v \rangle  =  \int_{{\R}^N}\left(\nabla u\nabla v+V(x)uv\right)\mathrm{d}x
     -\int_{{\R}^N}f(x, u) v\mathrm{d}x, \ \ \ \ \forall \ u, v\in E.
 \end{equation}
 Let
 \begin{equation}\label{Ph0}
   \Phi_0(u)=\frac{1}{2}\int_{{\R}^N}\left(|\nabla u|^2+V_0(x)u^2\right)\mathrm{d}x-\int_{{\R}^N}F_0(x, u)\mathrm{d}x,
      \ \ \ \ \forall \ u\in E,
 \end{equation}
 Then $\Phi_0$ is also of class $C^{1}(E, \R)$, and
 \begin{equation}\label{Ph0d}
   \langle \Phi_0'(u), v \rangle  =  \int_{{\R}^N}\left(\nabla u\nabla v+V_0(x)uv\right)\mathrm{d}x
     -\int_{{\R}^N}f_0(x, u) v\mathrm{d}x, \ \ \ \ \forall \ u, v\in E.
 \end{equation}
 In view of \x{Au} and \x{No0}, we have
 \begin{equation}\label{Ph1}
   \Phi_0(u)  =   \frac{1}{2}\left(\|u^{+}\|^2-\|u^{-}\|^2\right)-\int_{{\R}^N}F_0(x, u)\mathrm{d}x
 \end{equation}
 and
 \begin{equation}\label{Phd1}
   \langle \Phi_0'(u), u \rangle  =  \|u^{+}\|^2-\|u^{-}\|^2-\int_{{\R}^N}f_0(x, u) u\mathrm{d}x, \ \ \ \ \forall \ u=u^{-}+u^{+}\in E.
 \end{equation}

 \vskip4mm
 \par
    We set
 \begin{equation}\label{Psi}
   \Psi(u)=\int_{{\R}^N}[-V_1(x)u^2+F(x, u)]\mathrm{d}x, \ \ \ \ \forall \ u\in E.
 \end{equation}

 \vskip4mm
 \par
   Employing a standard argument, one can easily verify the following fact:

 \vskip4mm
 \par\noindent
 {\bf Lemma 2.2.}\ \ {\it Suppose that} (V1), (V2), (F1) {\it and} (F2) {\it  are satisfied. Then $\Psi$ is nonnegative,
 weakly sequentially lower semi-continuous, and $\Psi'$ is weakly sequentially continuous.}

 \vskip10mm

 {\section{The periodic case}}
 \setcounter{equation}{0}

 \vskip4mm
 \par
   In this section, we assume that $V$ and $f$ are 1-periodic in each of $x_1, x_2, \ldots, x_N$, i.e., (V) and
 (F3) are satisfied. In this case, $V_0=V$, $V_1=0$, $f_0=f$ and $f_1=0$. Thus, $\Phi_0(u)=\Phi(u)$.

 \vskip4mm
 \par\noindent
 {\bf Lemma 3.1.}(\cite[Lemma 2.4]{Ta4})\ \ {\it Suppose that} (V), (F1), (F2) {\it and} (WN) {\it are satisfied. Then}
 \begin{eqnarray}\label{L31}
   \Phi(u) & \ge & \Phi(tu+w)+\frac{1}{2}\|w\|^2+\frac{1-t^2}{2}\langle\Phi'(u), u \rangle-t\langle\Phi'(u), w \rangle,\nonumber\\
           &     & \ \ \ \ \ \ \ \ \ \ \ \ \ \ \ \forall \ u\in E, \ \ t\ge 0, \ \ w\in E^{-}.
 \end{eqnarray}

 \vskip4mm
 \par
   Define
 \begin{equation}\label{Ne-}
   \mathcal{N}^{-}  = \left\{u\in E\setminus E^{-} : \langle \Phi'(u), u \rangle=\langle \Phi'(u), v \rangle=0,
       \ \forall \ v\in E^{-} \right\}.
 \end{equation}
 The set $\mathcal{N}^{-}$ was first introduced by Pankov \cite{Pa}, which is a subset of the Nehari manifold
 $\mathcal{N}  = \left\{u\in E\setminus \{0\} : \langle \Phi'(u), u \rangle=0 \right\}.$

 \vskip4mm
 \par\noindent
 {\bf Corollary 3.2.}\ \ {\it Suppose that} (V), (F1), (F2) {\it and} (WN) {\it are satisfied. Then for
 $u\in \mathcal{N}^{-}$}
 \begin{equation}\label{2701}
   \Phi(u) \ge \Phi(tu+w)+\frac{1}{2}\|w\|^2, \ \ \ \ \forall \ t\ge 0, \ \ w\in E^{-}.
 \end{equation}

 \vskip4mm
 \par\noindent
 {\bf Corollary 3.3.}\ \ {\it Suppose that} (V), (F1), (F2) {\it and} (WN) {\it are satisfied. Then}
 \begin{eqnarray}\label{2801}
   \Phi(u) & \ge & \frac{t^2}{2}\|u\|^2-\int_{{\R}^N}F(x, tu^{+})\mathrm{d}x
                       +\frac{1-t^2}{2}\langle\Phi'(u), u \rangle+t^2\langle\Phi'(u), u^{-} \rangle,\nonumber\\
             &     & \ \ \ \ \ \ \ \ \ \ \ \ \ \ \ \ \ \ \ \ \ \ \ \ \ \ \ \ \ \ \ \ \ \ \ \ \forall \ u\in E, \ \ t\ge 0. \ \ \ \
 \end{eqnarray}

 \vskip4mm
 \par
    Analogous to the proof of \cite[Lemma 3.3]{Ta1}, it is easy to show the following lemma.

 \vskip4mm
 \par\noindent
 {\bf Lemma 3.4.}\ \ {\it Suppose that} (V), (F1), (F2) {\it and} (WS) {\it are satisfied.
 Then there exist a constant $c>0$ and a sequence $\{u_n\}\subset E$ satisfying}
 \begin{equation}\label{Ce0}
   \Phi(u_n)\rightarrow c, \ \ \ \ \|\Phi'(u_n)\|(1+\|u_n\|)\rightarrow 0.
 \end{equation}

 \vskip4mm
 \par\noindent
 {\bf Lemma 3.5.}\ \ {\it Suppose that} (V), (F1), (F2), (F3), (F4) {\it and} (WS) {\it are satisfied. Then any sequence
 $\{u_n\}\subset E$ satisfying
 \begin{equation}\label{Ce1}
   \Phi(u_n)\rightarrow c\ge 0, \ \ \ \ \langle\Phi'(u_n), u_n^{\pm} \rangle\rightarrow 0
 \end{equation}
 is bounded in $E$.}

 \vskip2mm
 \par\noindent
 {\bf Proof.} \ \ In view of \x{Ce1}, there exists a constant $C_2>0$ such that
 \begin{equation}\label{C2}
   C_2\ge \Phi(u_n)-\frac{1}{2}\langle\Phi'(u_n), u_n\rangle=\int_{\R^N}  \mathcal{F}(x, u_n)\mathrm{d}x.
 \end{equation}
 To prove the boundedness of $\{u_n\}$, arguing by contradiction, suppose that
 $\|u_n\| \to \infty$. Let $v_n=u_n/\|u_n\|$. Then $1=\|v_n\|^2$.
 If $\delta:=\limsup_{n\to\infty}\sup_{y\in \R^N}\int_{B_1(y)}|v_n^{+}|^2\mathrm{d}x=0$, then by Lions' concentration
 compactness principle \cite{Lio} or \cite[Lemma 1.21]{Wi}, $v_n^{+}\rightarrow 0$ in $L^{s}(\R^N)$ for $2<s<2^*$. Set
 $\kappa'=\kappa/(\kappa-1)$ and
 \begin{equation}\label{Omn}
   \Omega_n:=\left\{x\in \R^N : \frac{f(x, u_n)}{u_n}\le \bar{\Lambda}-\delta_0\right\}.
 \end{equation}
 Then using $\bar{\Lambda}\|v_n^{+}\|_2^2\le \|v_n^{+}\|^2$, one has
 \begin{eqnarray}\label{3903}
   \int_{\Omega_n}\frac{f(x, u_n)}{u_n}(v_n^{+})^2\mathrm{d}x\le (\bar{\Lambda}-\delta_0)\|v_n^{+}\|_2^2
   \le 1-\frac{\delta_0}{\bar{\Lambda}}.
 \end{eqnarray}
 On the other hand, by virtue of (F4), \x{C2} and the H\"older inequality, one can get that
 \begin{eqnarray}\label{3905}
  \int_{\R^N\setminus\Omega_n}\frac{f(x, u_n)}{u_n}(v_n^{+})^2\mathrm{d}x
   & \le & \left[\int_{\R^N\setminus\Omega_n}\left|\frac{f(x, u_n)}{u_n}\right|^{\kappa}\mathrm{d}x\right]^{1/\kappa}
             \|v_n^{+}\|_{2\kappa'}^2\nonumber\\
   & \le & C_2\left(\int_{\R^N\setminus\Omega_n} \mathcal{F}(x, u_n)\mathrm{d}x\right)^{1/\kappa}\|v_n^{+}\|_{2\kappa'}^2\nonumber\\
   & \le & C_3\|v_n^{+}\|_{2\kappa'}^2 =o(1).
 \end{eqnarray}
 $\mathcal{F}(x, u)\ge 0$ implies that $uf(x, u)\ge 0$. Hence, combining \x{3903} with \x{3905} and 
 making use of \x{Phd} and \x{Ce1}, we have
 \begin{eqnarray*}
  1+o(1)
   &  =  & \frac{\|u_n\|^2-\langle\Phi(u_n), u_n^{+}-u_n^{-}\rangle}{\|u_n\|^2}\\
   &  =  & \int_{u_n\ne 0}\frac{f(x, u_n)}{u_n}\left[(v_n^{+})^2-(v_n^{-})^2\right]\mathrm{d}x\\
   & \le & \int_{\Omega_n}\frac{f(x, u_n)}{u_n}(v_n^{+})^2\mathrm{d}x
            +\int_{\R^N\setminus\Omega_n}\frac{f(x, u_n)}{u_n}(v_n^{+})^2\mathrm{d}x\\
   & \le & 1-\frac{\delta_0}{\bar{\Lambda}}+o(1).
 \end{eqnarray*}
 This contradiction shows that $\delta>0$.

 \par
    Going if necessary to a subsequence, we may assume the existence of $k_n\in \Z^N$ such that
 $\int_{B_{1+\sqrt{N}}(k_n)}|v_n^{+}|^2\mathrm{d}x> \frac{\delta}{2}$. Let $w_n(x)=v_n(x+k_n)$. Since
 $V(x)$ is 1-periodic in each of $x_1, x_2, \ldots, $ $x_N$. Then
 $\|w_n\|=\|v_n\|=1$, and
 \begin{equation}\label{L62}
   \int_{B_{1+\sqrt{N}}(0)}|w_n^{+}|^2\mathrm{d}x> \frac{\delta}{2}.
 \end{equation}
 Passing to a subsequence, we have $w_n\rightharpoonup w$ in $E$, $w_n\rightarrow w$ in $L^{s}_{\mathrm{loc}}(\R^N)$,
 $2\le s<2^*$, $w_n\rightarrow w$ a.e. on $\R^N$. Obviously, \x{L62} implies that $w\ne 0$.

 \par
    Now we define $u_n^{k_n}(x)=u_n(x+k_n)$, then $u_n^{k_n}/\|u_n\|=w_n\rightarrow w$ a.e. on $\R^N$, $w\ne 0$.
 For $x\in \{y\in \R^N : w(y)\ne 0\}$, we have $\lim_{n\to\infty}|u_n^{k_n}(x)|=\infty$. Hence, it follows
 from \x{Ph}, \x{Ce1}, (F2), (F3), (WS) and Fatou's lemma that
 \begin{eqnarray*}
  0 &  =  & \lim_{n\to\infty}\frac{c+o(1)}{\|u_n\|^2} = \lim_{n\to\infty}\frac{\Phi(u_n)}{\|u_n\|^2}\\
    &  =  & \lim_{n\to\infty}\left[\frac{1}{2}\left(\|v_n^{+}\|^2-\|v_n^{-}\|^2\right)
             -\int_{\R^N}\frac{F(x, u_n^{k_n})}{(u_n^{k_n})^2}w_n^2\mathrm{d}x\right]\\
    & \le & \frac{1}{2}-\liminf_{n\to\infty}\int_{\R^N}\frac{F(x, u_n^{k_n})}{(u_n^{k_n})^2}w_n^2\mathrm{d}x
             \le \frac{1}{2}-\int_{\R^N}\liminf_{n\to\infty}\frac{F(x, u_n^{k_n})}{(u_n^{k_n})^2}w_n^2\mathrm{d}x  = -\infty.
 \end{eqnarray*}
 This contradiction shows that $\{u_n\}$ is bounded.
 \qed

 \vskip4mm
 \par\noindent
 {\bf Lemma 3.6.}(\cite[Theorem 1.2]{Ta5})\ \ {\it Assume that} (V), (F1), (F2), (F3), (WN) {\it and} (WS) {\it are satisfied. Then problem
 \x{SP} has a solution $u_0\in E$ such that $\Phi(u_0)=\inf_{\mathcal{N}^{-}}\Phi>0$.}

 \vskip4mm
 \par\noindent
 {\bf Proof of Theorem 1.3.} \  Combining Lemma 3.4 with Lemma 3.5, we can get that there exists a bounded sequence
 $\{u_n\}\subset E$ satisfying \x{Ce0}. Now the usual concentration-compactness argument suggests
 that $\Phi'(\bar{u})=0$ for some $\bar{u}\in E\setminus \{0\}$.
 \qed

 \vskip10mm

 {\section{The asymptotically periodic case}}
 \setcounter{equation}{0}

 \vskip4mm
 \par
    In this section, we always assume that $V$ satisfies (V1) and (V2).

 \vskip4mm
 \par\noindent
 {\bf Lemma 4.1.} \ \ {\it Suppose that} (V1), (V2), (F1) {\it and} (F2)  {\it are satisfied. Then
 there exists $\rho>0$ such that }
 \begin{equation}\label{m1}
    \hat{\kappa}:=\inf \left\{\Phi(u) : u\in E^{+}, \|u\|=\rho\right\}>0.
 \end{equation}

 \vskip2mm
 \par\noindent
 {\bf Proof.} \ \ Set $\Theta_0=\sup_{\R^N}[-V_1(x)]$. Let $\varepsilon_0=(\bar{\Lambda}-\Theta_0)/3$.
 Then (F1) and (F2) imply that there exists a constant $C_{\varepsilon_0}>0$ such that
 \begin{equation}\label{2501}
   F(x, t)\le \varepsilon_0|t|^2+C_{\varepsilon_0}|t|^{p}, \ \ \ \ \forall \ (x, t)\in \R^N\times \R.
 \end{equation}
 From \x{Ph}, \x{2501} and the Sobolev imbedding inequalities $\|u\|_p\le \gamma_p\|u\|$ and $\bar{\Lambda}\|u\|_2^2\le \|u\|^2$
 for $u\in E^{+}$, we have
 \begin{eqnarray}\label{2502}
   \Phi(u)
     &  =  & \frac{1}{2}\|u\|^2+\frac{1}{2}\int_{{\R}^N}V_1(x)u^2\mathrm{d}x
               -\int_{\R^N}F(x, u)\mathrm{d}x\nonumber\\
     & \ge & \frac{1}{2}\left[\|u\|^2-(\Theta_0+2\varepsilon_0)\|u\|_2^2\right]
               -C_{\varepsilon_0}\|u\|_p^p\nonumber\\
     & \ge & \frac{1}{2}\left(1-\frac{\Theta_0+2\varepsilon_0}{\bar{\Lambda}}\right)\|u\|^2
               -\gamma_p^pC_{\varepsilon_0}\|u\|^p, \ \ \ \ \forall \ u\in E^{+}.
 \end{eqnarray}
 This shows that there exists $\rho>0$ such that \x{m1} holds.

 \vskip4mm
 \par\noindent
 {\bf Lemma 4.2.}\ \ {\it Suppose that} (V1), (V2), (F1), (F2) {\it and} (WS) {\it are satisfied. Then for any
 $e\in E^{+}$, $\sup \Phi(E^{-}\oplus \R^{+} e)<\infty$, and there is $R_e>0$ such that}
 \begin{equation}\label{L321}
   \Phi(u) \le 0, \ \ \ \ \forall \   u\in E^{-}\oplus \R^{+} e, \ \ \|u\|\ge R_e.
 \end{equation}

 \vskip2mm
 \par\noindent
 {\bf Proof.} \ \ Arguing indirectly, provided that for some sequence $\{w_n+s_ne\}\subset
 E^{  -}\oplus \R e$ with $\|w_n+s_ne\|\rightarrow \infty$ such that $\Phi(w_n+s_ne) \ge 0$ for all $n\in \N$.
 Set $v_n=(w_n+s_ne)/\|w_n+s_ne\|=v_n^{  -}+t_ne$, then $\|v_n^{-}+t_ne\|=1$.
 Passing to a subsequence, we may assume that $t_n\rightarrow \bar{t}$, $v_n^{  -}\rightharpoonup v^{  -}$,
 and $v_n^{-}\rightarrow v^{-}$ a.e. on $\R^N$. Hence, it follows from (V2) and \x{Ph} that
 \begin{eqnarray}\label{L402}
   0 & \le & \frac{\Phi(w_n+s_ne)}{\|w_n+s_nu\|^2}\nonumber\\
     &  =  & \frac{t_n^2}{2}\|e\|^2-\frac{1}{2}\|v_n^{-}\|^2+\frac{1}{2}\int_{{\R}^N}V_1(x)(v_n^{-}+t_ne)^2\mathrm{d}x
              -\int_{\R^N}\frac{F(x, w_n+s_ne)}{\|w_n+s_ne\|^2}\mathrm{d}x\nonumber\\
     & \le & \frac{t_n^2}{2}\|e\|^2-\frac{1}{2}\|v_n^{-}\|^2
              -\int_{\R^N}\frac{F(x, w_n+s_ne)}{\|w_n+s_ne\|^2}\mathrm{d}x.
 \end{eqnarray}

 \par
   If $\bar{t}=0$, then it follows from \x{L402} that
 $$
   0 \le \frac{1}{2}\|v_n^{  -}\|^2+\int_{\R^N}\frac{F(x, w_n+s_ne)}{\|w_n+s_ne\|^2}\mathrm{d}x
     \le \frac{t_n^2}{2}\|e\|^2 \rightarrow 0,
 $$
 which yields $\|v_n^{-}\|\rightarrow 0$, and so $1=\|v_n^{-}+t_ne\|^2 \rightarrow 0$, a contradiction.

 \par
   If $\bar{t}\ne 0$, then it follows from \x{L402} and (WS) that
 \begin{eqnarray*}
  0 & \le & \limsup_{n\to\infty}\left[\frac{t_n^2}{2}\|e\|^2-\frac{1}{2}\|v_n^{-}\|^2
              -\int_{\R^N}\frac{F(x, w_n+s_ne)}{\|w_n+s_ne\|^2}\mathrm{d}x\right]\\
    & \le & \frac{\bar{t}^2}{2}\|e\|^2
             -\liminf_{n\to\infty}\int_{\R^N}\frac{F(x, w_n+s_ne)}{(w_n+s_ne)^2}\left(v_n^{-}+t_ne\right)^2\mathrm{d}x\\
    & \le & \frac{\bar{t}^2}{2}\|e\|^2
             -\int_{\R^N}\liminf_{n\to\infty}\frac{F(x, w_n+s_ne)}{(w_n+s_ne)^2}\left(v_n^{-}+t_ne\right)^2\mathrm{d}x\\
    &  =  & -\infty,
 \end{eqnarray*}
 a contradiction.
 \qed

 \vskip4mm
 \par\noindent
 {\bf Corollary 4.3.}\ \ {\it Suppose that} (V1), (V2), (F1), (F2) {\it and} (WS) {\it are satisfied.
 Let $e\in E^{+}$ with $\|e\|=1$. Then there is a $r>\rho$ such that $\sup \Phi(\partial Q)\le 0$, where}
 \begin{equation}\label{Q1}
   Q=\left\{w+se : w\in E^{-}, \ s\ge 0,\ \|w+se\|\le r\right\}.
 \end{equation}

 \vskip4mm
 \par
   Let $m_0:=\inf_{\mathcal{N}^{0}}\Phi_0$, where
 \begin{equation}\label{Ne0}
   \mathcal{N}^{0}  = \left\{u\in E\setminus E^{-} : \langle \Phi_0'(u), u \rangle=\langle \Phi_0'(u), v \rangle=0,
       \ \forall \ v\in E^{-} \right\}.
 \end{equation}

 \vskip4mm
 \par\noindent
 {\bf Lemma 4.4.}\ \ {\it Suppose that} (V1), (V2), (F1), (F2) {\it and} (F5$'$)  {\it are satisfied.
 Then there exist a constant $c_*\in [\hat{\kappa}, m_0)$ and a sequence $\{u_n\}\subset E$ satisfying}
 \begin{equation}\label{Ce2}
   \Phi(u_n)\rightarrow c_*, \ \ \ \ \|\Phi'(u_n)\|(1+\|u_n\|)\rightarrow 0,
 \end{equation}

 \vskip2mm
 \par\noindent
 {\bf Proof.} \ \ Employing Lemma 3.6, there exists a $\bar{u}\in E$ such that $\bar{u}\ne 0$ on $B_{1+\sqrt{N}}(0)$,
 $\Phi_0'(\bar{u})=0$ and $\Phi_0(\bar{u})=m_0$. Set $\hat{E}(\bar{u})=E^{-}\oplus \R^{+} \bar{u}$ and $\zeta_0:=\sup\{\Phi(u) : u\in \hat{E}(\bar{u})\}$. Lemma 4.1 implies that $\zeta_0\ge \hat{\kappa}>0$. By (V2), (F5$'$), \x{Ph}, \x{Ph0} and Corollary 3.2, we have
 \begin{equation}\label{2601}
   \Phi(u)\le \Phi_0(u)\le m_0, \ \ \ \ \forall \ u\in \hat{E}(\bar{u}).
 \end{equation}
 Hence, $\zeta_0\le m_0$. If $\zeta_0 = m_0$, then there is a sequence $\{u_n\}$ with $u_n=w_n+s_n\bar{u}\in \hat{E}(\bar{u})$
 such that
 \begin{equation}\label{2603}
   m_0-\frac{1}{n}<\Phi(u_n)=\Phi(w_n+s_n\bar{u})\le m_0.
 \end{equation}
 It follows from Lemma 4.2 and \x{2603} that $\{s_n\}\subset \R$ and $\{w_n\}\subset E^{-}$ are bounded.
 Passing to a subsequence, we have $s_n\rightarrow \bar{s}$ and $w_n\rightharpoonup \bar{w}$ in $E$. It
 is easy to see that $\bar{s}>0$. It follows from \x{Ph}, \x{Ph0} and Corollary 3.2 that
 \begin{eqnarray*}
   m_0-\frac{1}{n}
    &  <  & \Phi(u_n)=\Phi_0(u_n)+\frac{1}{2}\int_{\R^N}V_1(x)u_n^2\mathrm{d}x-\int_{\R^N}F_1(x, u_n)\mathrm{d}x\\
    & \le & m_0-\frac{1}{2}\|w_n\|^2-\frac{1}{2}\int_{\R^N}[-V_1(x)u_n^2+2F_1(x, u_n)]\mathrm{d}x,
 \end{eqnarray*}
 which yields that
 \begin{equation}\label{2605}
   \frac{1}{2}\|w_n\|^2+\frac{1}{2}\int_{\R^N}[-V_1(x)(w_n+s_n\bar{u})^2+2F_1(x, w_n+s_n\bar{u})]\mathrm{d}x \le \frac{1}{n}.
 \end{equation}
 According to Fatou's lemma and the weakly lower semi-continuous of the norm, one gets that
 \begin{equation}\label{2607}
   \frac{1}{2}\|\bar{w}\|^2+\frac{1}{2}\int_{\R^N}[-V_1(x)(\bar{w}+\bar{s}\bar{u})^2+2F_1(x, \bar{w}+\bar{s}\bar{u})]\mathrm{d}x=0.
 \end{equation}
 This, together with (F5$'$), implies that $\bar{w}=0$ and $\bar{u}= 0$ on $B_{1+\sqrt{N}}(0)$, a contradiction.
 Therefore, $\zeta_0 \in [\hat{\kappa}, m_0)$. In view of Lemmas 2.1, 2.2, 4.1 and Corollary 4.3, there exist a constant
 $c_*\in [\hat{\kappa}, m_0)$ and a sequence $\{u_n\}\subset E$ satisfying \x{Ce2}.
 \qed

 \vskip4mm
 \par\noindent
 {\bf Lemma 4.5.}\ \ {\it Suppose that} (V1), (V2), (F1), (F2), (F4) {\it and} (F5$'$) {\it are satisfied. Then any sequence
 $\{u_n\}\subset E$ satisfying \x{Ce1} is bounded in $E$.}

 \vskip2mm
 \par\noindent
 {\bf Proof.} \ \ Given the condition \x{Ce1}, there exists a constant $C_2>0$ such that \x{C2} holds.
 To prove the boundedness of $\{u_n\}$, arguing by contradiction, suppose that
 $\|u_n\| \to \infty$. Let $v_n=u_n/\|u_n\|$. Then $1=\|v_n\|^2$. Passing to a subsequence, we have
 $v_n\rightharpoonup \bar{v}$ in $E$. There are two possible cases: i). $\bar{v}= 0$ and ii). $\bar{v}\ne 0$.

 \vskip2mm
 \par
   Case 1). $\bar{v}=0$, i.e. $v_n\rightharpoonup 0$ in $E$. Then $v_n^{+}\rightarrow 0$ and
 $v_n^{-}\rightarrow 0$ in $L^{s}_{\mathrm{loc}}(\R^N)$, $2\le s<2^*$
 and $v_n^{+}\rightarrow 0$ and $v_n^{-}\rightarrow 0$ a.e. on $\R^N$. From (V1), it is easy
 to show that
 \begin{equation}\label{4501}
   \lim_{n\to\infty}\int_{\R^N}V_1(x)(v_n^{+})^2\mathrm{d}x=\lim_{n\to\infty}\int_{\R^N}V_1(x)(v_n^{-})^2\mathrm{d}x=0.
 \end{equation}

 \par
    If $\delta:=\limsup_{n\to\infty}\sup_{y\in \R^N}\int_{B_1(y)}|v_n^{+}|^2\mathrm{d}x=0$, then by Lions' concentration
 compactness principle \cite{Lio} or \cite[Lemma 1.21]{Wi}, $v_n^{+}\rightarrow 0$ in $L^{s}(\R^N)$ for $2<s<2^*$. Set
 $\Omega_n$ as \x{Omn}. Then \x{3903} and \x{3905} hold also. Combining \x{3903} with \x{3905} and using \x{Phd}, \x{Ce1}
 and \x{4501}, we have
 \begin{eqnarray*}
  1+o(1)
   &  =  & \frac{\|u_n\|^2-\langle\Phi(u_n), u_n^{+}-u_n^{-}\rangle}{\|u_n\|^2}\\
   &  =  & -\int_{\R^N}V_1(x)\left[(v_n^{+})^2-(v_n^{-})^2\right]\mathrm{d}x
            +\int_{u_n\ne 0}\frac{f(x, u_n)}{u_n}\left[(v_n^{+})^2-(v_n^{-})^2\right]\mathrm{d}x\\
   & \le & -\int_{\R^N}V_1(x)(v_n^{+})^2\mathrm{d}x
            +\int_{\Omega_n}\frac{f(x, u_n)}{u_n}(v_n^{+})^2\mathrm{d}x
            +\int_{\R^N\setminus\Omega_n}\frac{f(x, u_n)}{u_n}(v_n^{+})^2\mathrm{d}x\\
   & \le & 1-\frac{\delta_0}{\bar{\Lambda}}+o(1).
 \end{eqnarray*}
 This contradiction shows that $\delta>0$.

 \par
    Going if necessary to a subsequence, we may assume the existence of $k_n\in \Z^N$ such that
 $\int_{B_{1+\sqrt{N}}(k_n)}|v_n^{+}|^2\mathrm{d}x> \frac{\delta}{2}$. Let $w_n(x)=v_n(x+k_n)$. Since
 $V_0(x)$ is 1-periodic in each of $x_1, x_2, \ldots, $ $x_N$. Then
 \begin{equation}\label{l13}
   \int_{B_{1+\sqrt{N}}(0)}|w_n^{  +}|^2\mathrm{d}x> \frac{\delta}{2}.
 \end{equation}
 Now we define $\tilde{u}_n(x)=u_n(x+k_n)$, then $\tilde{u}_n/\|u_n\|=w_n$ and $\|w_n\|=\|v_n\|=1$.
 Passing to a subsequence, we have $w_n\rightharpoonup w$ in $E$, $w_n\rightarrow w$ in $L^{s}_{\mathrm{loc}}(\R^N)$,
 $2\le s<2^*$ and $w_n\rightarrow w$ a.e. on $\R^N$. Obviously, \x{l13} implies that $w\ne 0$. Hence, it follows from \x{Ce1},
 (F5$'$) and Fatou's lemma that
 \begin{eqnarray*}
  0 &  =  & \lim_{n\to\infty}\frac{c+o(1)}{\|u_n\|^2} = \lim_{n\to\infty}\frac{\Phi(u_n)}{\|u_n\|^2}\nonumber\\
    &  =  & \lim_{n\to\infty}\left[\frac{1}{2}\left(\|v_n^{+}\|^2-\|v_n^{-}\|^2\right)
              +\frac{1}{2}\int_{{\R}^N}V_1(x)v_n^2\mathrm{d}x -\int_{\R^N}\frac{F(x+k_n,\tilde{u}_n)}{\tilde{u}_n^2}w_n^2dx\right]\nonumber\\
    & \le & \frac{1}{2}-\liminf_{n\to\infty}\int_{\R^N}\frac{F_0(x, \tilde{u}_n)}{\tilde{u}_n^2}w_n^2dx
             \le \frac{1}{2}-\int_{\R^N}\liminf_{n\to\infty}\frac{F_0(x, \tilde{u}_n)}{\tilde{u}_n^2}w_n^2dx\nonumber\\
    &  =  & -\infty,
 \end{eqnarray*}
 which is a contradiction.

 \vskip2mm
 \par
   Case ii). $\bar{v}\ne 0$. In this case, we can also deduce a contradiction by a standard argument.

 \vskip2mm
 \par
 Cases i) and ii) both show that $\{u_n\}$ is bounded in $E$.
 \qed

 \vskip4mm
 \par\noindent
 {\bf Proof of Theorem 1.5.} \ \ It is easy to see that (F5$'$) implies (WS). Applying Lemmas 4.4 and 4.5,
 we obtain that there exists a bounded sequence $\{u_n\}\subset E$ satisfying \x{Ce2}. Passing to a subsequence,
 we have $u_n\rightharpoonup \bar{u}$ in $E$. Next, we prove $\bar{u}\ne 0$.

 \par
   Arguing by contradiction, suppose that $\bar{u}=0$, i.e. $u_n\rightharpoonup 0$ in $E$, and so $u_n\rightarrow 0$
 in $L^{s}_{\mathrm{loc}}(\R^N)$, $2\le s<2^*$ and $u_n\rightarrow 0$ a.e. on $\R^N$. By (V1) and (F5$'$),
 it is easy to show that
 \begin{equation}\label{T32}
   \lim_{n\to\infty}\int_{\R^N}V_1(x)u_n^2dx=0, \ \ \ \
   \lim_{n\to\infty}\int_{\R^N}V_1(x)u_nvdx=0, \ \ \ \ \forall \ v\in E
 \end{equation}
 and
 \begin{equation}\label{T34}
   \lim_{n\to\infty}\int_{\R^N}F_1(x, u_n)dx=0, \ \ \ \ \lim_{n\to\infty}\int_{\R^N}f_1(x, u_n)vdx=0,
      \ \ \ \ \forall \ v\in E.
 \end{equation}
 Note that
 \begin{equation}\label{T35}
   \Phi_0(u)=\Phi(u)-\frac{1}{2}\int_{\R^N}V_1(x)u^2dx+\int_{\R^N}F_1(x, u)dx, \ \ \ \ \forall \ u\in E
 \end{equation}
 and
 \begin{equation}\label{T36}
   \langle\Phi_0'(u), v\rangle=\langle\Phi'(u), v\rangle-\int_{\R^N}V_1(x)uvdx+\int_{\R^N}f_1(x, u)vdx, \ \ \ \ \forall \ u, v\in E.
 \end{equation}
 From \x{Ce2}, \x{T32}-\x{T36}, one can get that
 \begin{equation}\label{T37}
   \Phi_0(u_n)\rightarrow c_*, \ \ \ \ \|\Phi_0'(u_n)\|(1+\|u_n\|)\rightarrow 0.
 \end{equation}

 \par
   A standard argument shows that $\{u_n\}$ is non-vanishing sequence.
 Going if necessary to a subsequence, we may assume the existence of $k_n\in \Z^N$ such that
 $\int_{B_{1+\sqrt{N}}(k_n)}|u_n|^2\mathrm{d}x> \frac{\delta}{2}$ for some $\delta>0$. Let $v_n(x)=u_n(x+k_n)$.
 Then
 \begin{equation}\label{L74}
   \int_{B_{1+\sqrt{N}}(0)}|v_n|^2\mathrm{d}x> \frac{\delta}{2}.
 \end{equation}
 Since $V_0(x)$ and $f_0(x, u)$ are periodic on $x$, we have $\|v_n\|=\|u_n\|$ and
 \begin{equation}\label{L75}
   \Phi_0(v_n)\rightarrow c_*, \ \ \ \ \|\Phi_0'(v_n)\|(1+\|v_n\|)\rightarrow 0.
 \end{equation}
 Passing to a subsequence, we have $v_n\rightharpoonup \bar{v}$ in $E$, $v_n\rightarrow \bar{v}$ in $L^{s}_{\mathrm{loc}}(\R^N)$,
 $2\le s<2^*$ and $v_n\rightarrow \bar{v}$ a.e. on $\R^N$. Obviously, \x{L74} and \x{L75} imply that $\bar{v}\ne 0$
 and $\Phi_0'(\bar{v})=0$. This shows that $\bar{v}\in \mathcal{N}^{0}$  and so $\Phi_0(\bar{v})\ge m_0$.
 On the other hand, by using \x{L75}, (WN) and Fatou's lemma, we have
 \begin{eqnarray*}
   m_0 &  >  & c_*=\lim_{n\to\infty}\left[\Phi_0(v_n)-\frac{1}{2}\langle\Phi_0'(v_n), v_n \rangle\right]
                = \lim_{n\to\infty}\int_{\R^N}\left[\frac{1}{2}f_0(x, v_n)v_n-F_0(x, v_n)\right]\mathrm{d}x\\
       & \ge & \int_{\R^N}\lim_{n\to\infty}\left[\frac{1}{2}f_0(x, v_n)v_n-F_0(x, v_n)\right]\mathrm{d}x
                =  \int_{\R^N}\left[\frac{1}{2}f_0(x, \bar{v})\bar{v}-F_0(x, \bar{v})\right]\mathrm{d}x\\
       &  =  & \Phi_0(\bar{v})-\frac{1}{2}\langle\Phi_0'(\bar{v}), \bar{v} \rangle = \Phi_0(\bar{v})\ge m_0.
 \end{eqnarray*}
 This contradiction implies that $\bar{u}\ne 0$. It is obvious that that $\bar{u}\in E$ is a nontrivial solution for
 problem \x{SP}.
 \qed

 {}
\end{document}